\documentclass[a4paper,10pt]{article}
\usepackage[utf8]{inputenc} 
\usepackage{amssymb,latexsym}
\usepackage{amsmath,amsthm}
\usepackage{graphicx}
\usepackage{color,epsfig}
\usepackage{physics}
\usepackage{hyperref}

\usepackage{amsmath,amssymb} 
\usepackage{graphicx}
\usepackage{hyperref}      
\usepackage{xcolor}
\usepackage{colortbl}
\usepackage{blindtext}
\usepackage{array}
\usepackage{amsmath}
\usepackage{amsthm}
\usepackage{amssymb}
\usepackage{tikz}
\usepackage{enumitem}
\usepackage{mathrsfs}
\usepackage{wrapfig}

\usepackage[utf8]{inputenc} 

\usepackage{mathtools}
\input epsf
\usepackage{epsfig}
\usepackage{graphicx}
\usepackage{epstopdf}
\usepackage{subcaption}
\usepackage[overload]{empheq}
\usepackage{multicol}
\usepackage{currvita}
\usepackage{tabu}
\usepackage{float,graphicx}
\usepackage[justification=centering]{caption}
\usepackage{caption}
\usepackage{amsmath}
\usepackage[utf8]{inputenc}
\usepackage{amssymb}
\usepackage{physics}
\usepackage{changepage}
\usepackage{relsize}
\usepackage{svg}
\usepackage{amsmath}
\usepackage{etoolbox,refcount}
\usepackage{multicol}

\newtheorem{theorem}{Theorem}
\newtheorem{lemma}{Lemma}
\newtheorem{corollary}{Corollary}

\def\qed{\ifhmode\unskip\nobreak\fi\quad\ifmmode\Box\else$\Box$\fi}

\newcounter{countitems}
\newcounter{nextitemizecount}
\newcommand{\setupcountitems}{%
  \stepcounter{nextitemizecount}%
  \setcounter{countitems}{0}%
  \preto\item{\stepcounter{countitems}}%
}
\makeatletter
\newcommand{\computecountitems}{%
  \edef\@currentlabel{\number\c@countitems}%
  \label{countitems@\number\numexpr\value{nextitemizecount}-1\relax}%
}
\newcommand{\nextitemizecount}{%
  \getrefnumber{countitems@\number\c@nextitemizecount}%
}
\newcommand{\previtemizecount}{%
  \getrefnumber{countitems@\number\numexpr\value{nextitemizecount}-1\relax}%
}
\makeatother    
{\end{itemize}%
\unskip\computecountitems\ifnumcomp{\previtemizecount}{>}{3}{\end{multicols}}{}}

\title{Note on the chromatic number of Minkowski planes:\\
the regular polygon case}

\author{\sl Panna Gehér
\thanks{ Eötvös Loránd University, Budapest; \texttt{geherpanni@student.elte.hu}.}
}

\date{}

\begin{document}

\maketitle

\begin{abstract}
The famous Hadwiger-Nelson problem asks for the minimum number of colors needed to color the points of the Euclidean plane so that no two points unit distance apart are assigned the same color.
In this note we consider a variant of the problem in Minkowski metric planes, where the unit circle is a regular polygon of even and at most $22$ vertices.
We present a simple lattice–sublattice coloring scheme that uses $6$ colors, proving that the chromatic number of the Minkowski planes above are at most 6. This result is new for regular polygons having more than 8 vertices.
\end{abstract}

\section{Introduction}
In 1950, Nelson raised the following question:
What is the minimum number of colors that are needed to color the Euclidean plane so that no two points of the same color determine unit distance? 
We refer to such a coloring with $k$ color classes as a proper $k$-coloring.
Thus Nelson's question asks for the smallest $k$ value such that the plane can be properly $k$-colored. This value is known as the chromatic number of the Euclidean plane, and is denoted by $\chi(\mathbb{R}^2)$. Immediately after the question was raised the following easy-to-get bounds were established:
 $$4 \leq \chi(\mathbb{R}^2) \leq 7.$$

\begin{figure}[htp!]
\centering
\includegraphics[width=11.3cm]{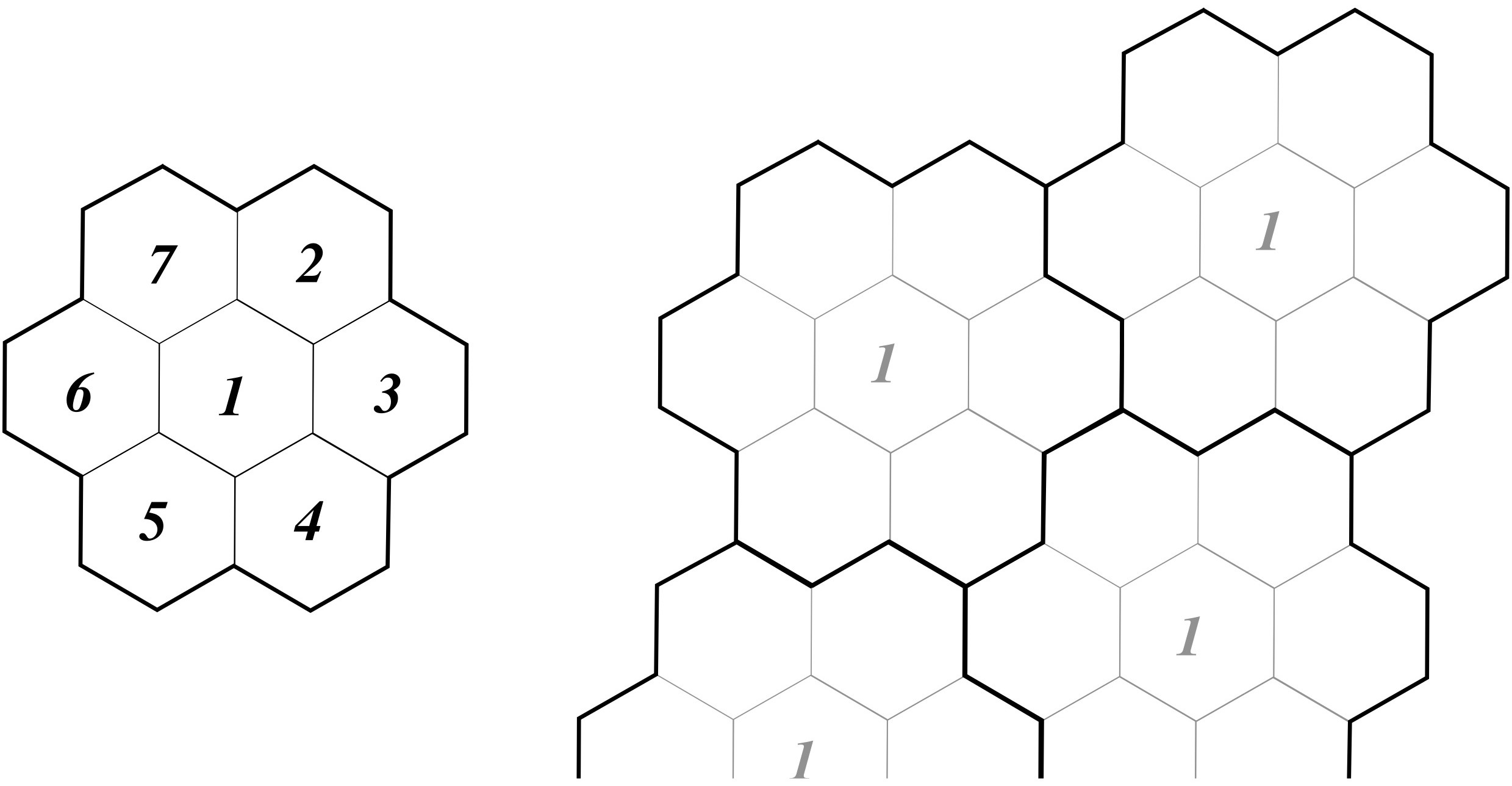}
\caption{\small{A proper $7$-coloring of the Euclidean plane shown by Isbell}}
\label{Isbell_kep}
\end{figure}

The lower bound is due to Moser \cite{Moser} who constructed a unit-distance graph (that is a graph whose edges connect vertices unit distance apart) with chromatic number $4$. The upper bound is due to Isbell \cite{Isbell} who considered a tilting of the plane by translates of a regular hexagon with diameter slightly less than one and defined a periodic proper $7$-coloring shown in Figure \ref{Isbell_kep}.

Despite the numerous attempts to improve these bounds, only little progress was made for more that $60$ years -- for a historical survey on the problem see \cite{Soifer}. However, in 2018 mathematicians were shaken when biologist de Grey \cite{de-Grey} constructed a $5$-chromatic unit-distance graph, proving that the chromatic number of the plane is at least $5$. Shortly afterwards Exoo and Ismailescu \cite{Exoo-Ismailescu} independently published another proof.\\
The problem has regained a lot of attention since the breakthrough and a Polymath project was launched with the main goal of creating a human-verifiable proof of the new result. Although the proofs are still relying on computers, quite some progress has been made: while the distance graph published by de Grey had a total of $1581$ vertices, the current known smallest example consists only $509$ \cite{Parts}.\\

As a consequence of the breakthrough many variations of the problem have gained more attention in the last couple of years.
One interesting research area is generalizing the question to Minkowski planes:
Let $C$ be a two-dimensional centrally symmetric bounded convex domain centered at the origin and let $(\mathbb{R}^2, \text{ } C)$ denote the Minkowski plane where $C$ determines the unit circle; the $C$-norm of an $x \in \mathbb{R}^2$ point is the value:
$$\norm{x}_C := \min \left\{ \lambda \in \mathbb{R}^{+}| x\in \lambda C\right\} .$$
The $C$-distance of two points $x$ and $y$ is defined by $\norm{x - y}_C$. Naturally, the chromatic number of the Minkowski planes -- denoted by $\chi(\mathbb{R}^2,C)$ -- is the minimum number of colors needed to color the points of $\mathbb{R}^2$ such that no monochromatic point pair determines a unit $C$-distance. 
The main result concerning the chromatic number of Minkowski planes is due to Chilakamarri \cite{Chilakamarri}: by extending the arguments of Moser and the construction of Isbell he proved that the bounds $$4 \leq \chi(\mathbb{R}^2, C) \leq 7$$ hold for all centrally symmetric bounded convex domain $C$.\\

\begin{figure}[htp]
\centering
\includegraphics[width=7.2cm]{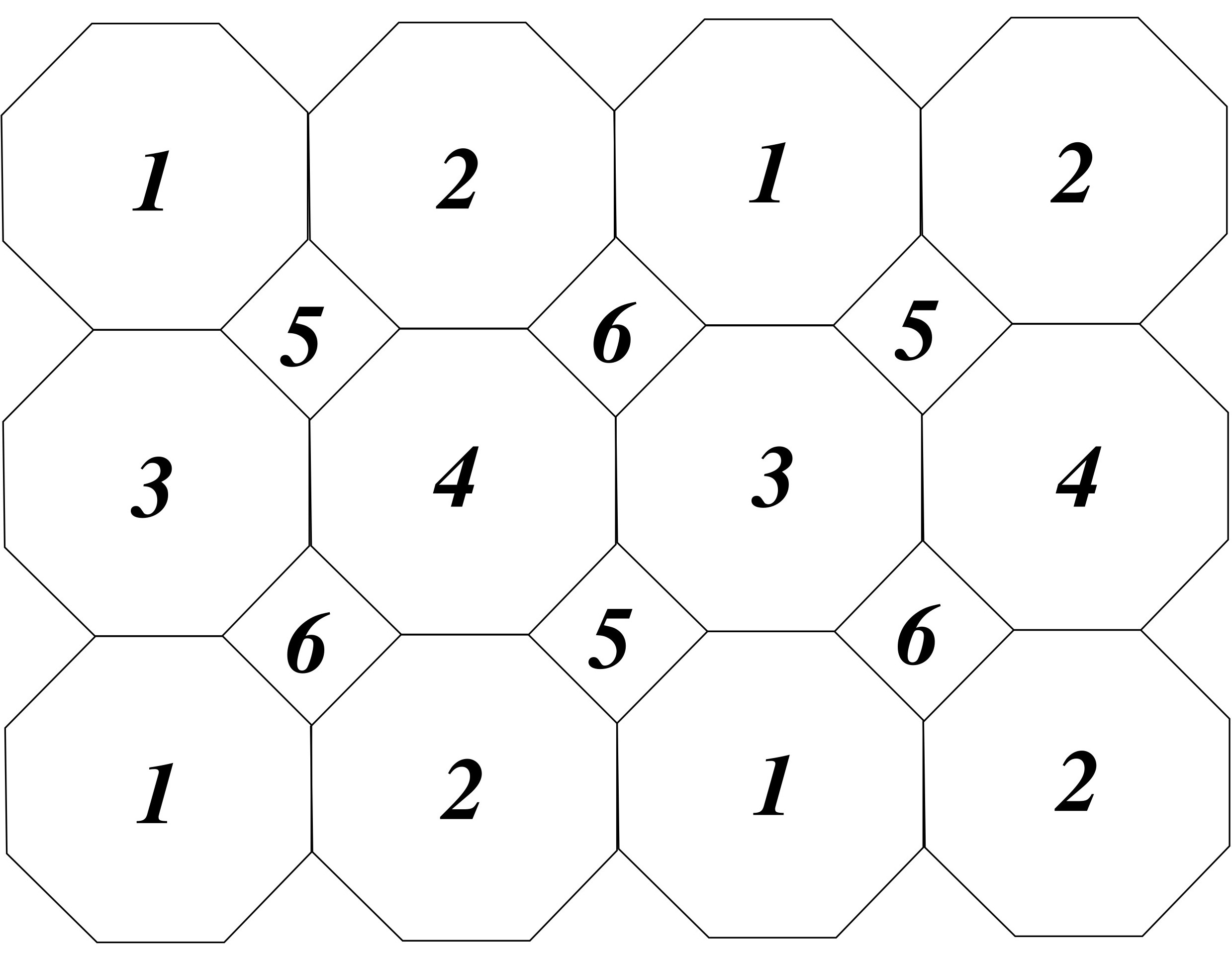}
\caption{\small{Chilakamarri's proper $6$-coloring of the Minkowski plane equipped with the regular octagon metric}}
\label{Chilakamarri_6szin}
\end{figure}

An interesting special case of the above problem is when the unit circle is a regular polygon of even number of vertices. Or more generally we can consider any affine image of a regular polygon since the problem itself is affine invariant.
The study of the chromatic number of such normed planes was also initiated by
Chilakamarri who considered the cases of regular polygons with few vertices. He gave a tile-based proper $4$-coloring for the parallelogram and the symmetric hexagon's case, proving that the answer here is exactly $4$. 
He also gave a proper $6$-coloring in case $C$ is a centrally symmetric octagon: he considered a packing of $C/2$, that is a packing of circles of radius half. He showed that the translates of $C/2$ can be colored using $4$ colors and two more colors can take care of the remaining squares (for details see Figure \ref{Chilakamarri_6szin}).
It is worth noting that we have to be careful with choosing the colors of boundary points of the octagons as
no antipodal point pair can have the same color.

Chilakamarri asked whether or not the chromatic number of a plane equipped with the norm defined by the regular octagon or any centrally symmetric octagon is exactly $4$.
No progress was made until very recently when Exoo, Fisher and Ismailescu \cite{Exoo-Fisher-Ismailescu} answered his question negatively:
they constructed $5$-chromatic unit-distance graphs in the cases of regular polygons with $8$, $10$, and $12$ vertices.
Together with the Euclidean case these are the only known examples of  a normed plane with chromatic number at least $5$.
Table \ref{table} summarizes the mentioned results:
\begin{table}[h]
  \begin{tabular}{  m{9.5cm} | m{2.5cm}   }
\hspace{0.2cm} Unit circle $C$ &  \hspace{0.5cm}  $\chi(\mathbb{R}^2,  C)$ \\ \hline
 
    \begin{itemize}[leftmargin=0.35cm, label={}]
        \item Parallelogram, centrally symmetric hexagon (see \cite{Chilakamarri})
        \vspace{-0.25cm}
        \item Regular octagon (see \cite{Exoo-Fisher-Ismailescu}, and \cite{Chilakamarri})
         \vspace{-0.25cm}
        \item Regular decagon, regular dodecagon (see \cite{Exoo-Fisher-Ismailescu} and \cite{Chilakamarri})
         \vspace{-0.25cm}
        \item Euclidean circle (see \cite{de-Grey}, \cite{Exoo-Ismailescu} and \cite{Isbell})
        \vspace{-0.25cm}
        \item Arbitrary symmetric convex domain (see \cite{Chilakamarri})
    \end{itemize}
    &    
    \begin{itemize} [leftmargin=0.7cm, label={}]
        \item $4$
          \vspace{-0.25cm}
        \item $5$ or $6$     
         \vspace{-0.25cm}
        \item $5$, $6$ or $7$ 
        \vspace{-0.25cm}
        \item $5$, $6$ or $7$ 
        \vspace{-0.25cm}
        \item $4$, $5$, $6$ or $7$ 

    \end{itemize}
   
  \end{tabular}

\caption{Possible values of the chromatic numbers of Minkowski planes}
\label{table}
\end{table}

In this note we extend Chilakamarri's result for regular octagons to regular polygons with at most $22$ vertices by giving a simple lattice–sublattice coloring scheme that uses only $6$ colors. It also slightly strengthens the result of Chilakamarri as our colorings are regular: 
We call a proper $k$-coloring of $\mathbb{R}^n$ with color classes $C_1$, $C_2$ \dots $C_k$ \textit{regular}, if there exist vectors $v_1 \dots$, $v_k$ such that $C_i=C_1+v_i$ for $i=1 \dots k$, that is the color classes are translates of each other.
Now we state our main theorem:

\begin{theorem} \label{theorem1}
Let $C$ be a regular polygon with an even number of vertices. In case $C$ has at most $22$ vertices then there exists a regular proper $6$-coloring of the Minkowski plane equipped with the $C$-metric such that no points unit $C$-distance apart are identically colored. Hence, $$\chi(\mathbb{R}^2, C) \leq 6$$ holds for any regular $2k$-gon $C$ where $k\leq 11$.
\end{theorem}

It follows that in the case of a regular decagon and dodecagon -- similarly to the octagon's case -- the chromatic number is either $5$ or $6$.

\newpage

In Section \ref{coloring} we describe the coloring scheme used in all cases and give some details and figures about the proofs. As an example, in Section \ref{proof} all computations are given in the regular dodecagon's case and calculations for the $22$-gon can be found in the Appendix.
Finally, in Section \ref{Szlam} we describe a consequence of Theorem \ref{theorem1} that considers a closely related asymmetric Ramsey-type question raised by Szlam \cite{Szlam}.


\section{The coloring scheme} \label{coloring}

\begin{figure}[htp!]
\centering
\includegraphics[width=4.5cm]{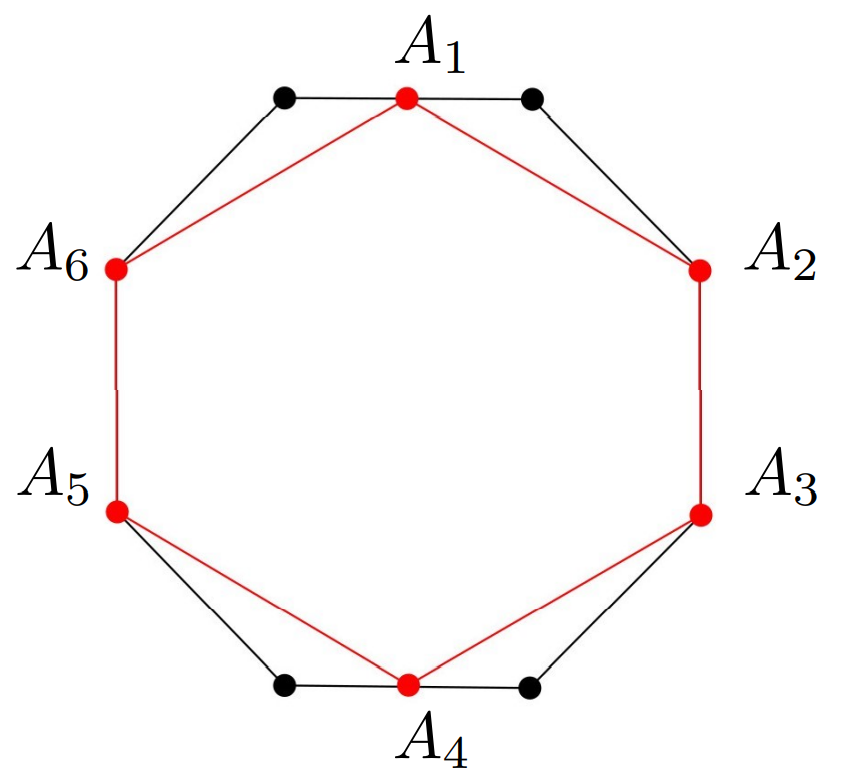}
\vspace{0.07cm}
\caption{\small{The centrally symmetric hexagons inscribed in $C/2$ in the case of the regular decagon}}
\label{hatszog_8}
\end{figure}

Let $C$ be a regular octagon first and define a symmetric convex hexagon $H$ inscribed in $C/2$ as follows: Choose two opposite sides of $C/2$ and form a hexagon using their four endpoints and two additional boundary points of $C/2$. The choice of the additional points can be made in various ways, here we simply chose the ones that halve the boundary line of $C/2$ connecting the chosen sides. Denote the vertices of $H$ by $A_i$ ($i=1 \dots 6$) in a clockwise order as shown in Figure \ref{hatszog_8}.

\begin{figure}[htp!]
\centering
\includegraphics[width=12.8cm]{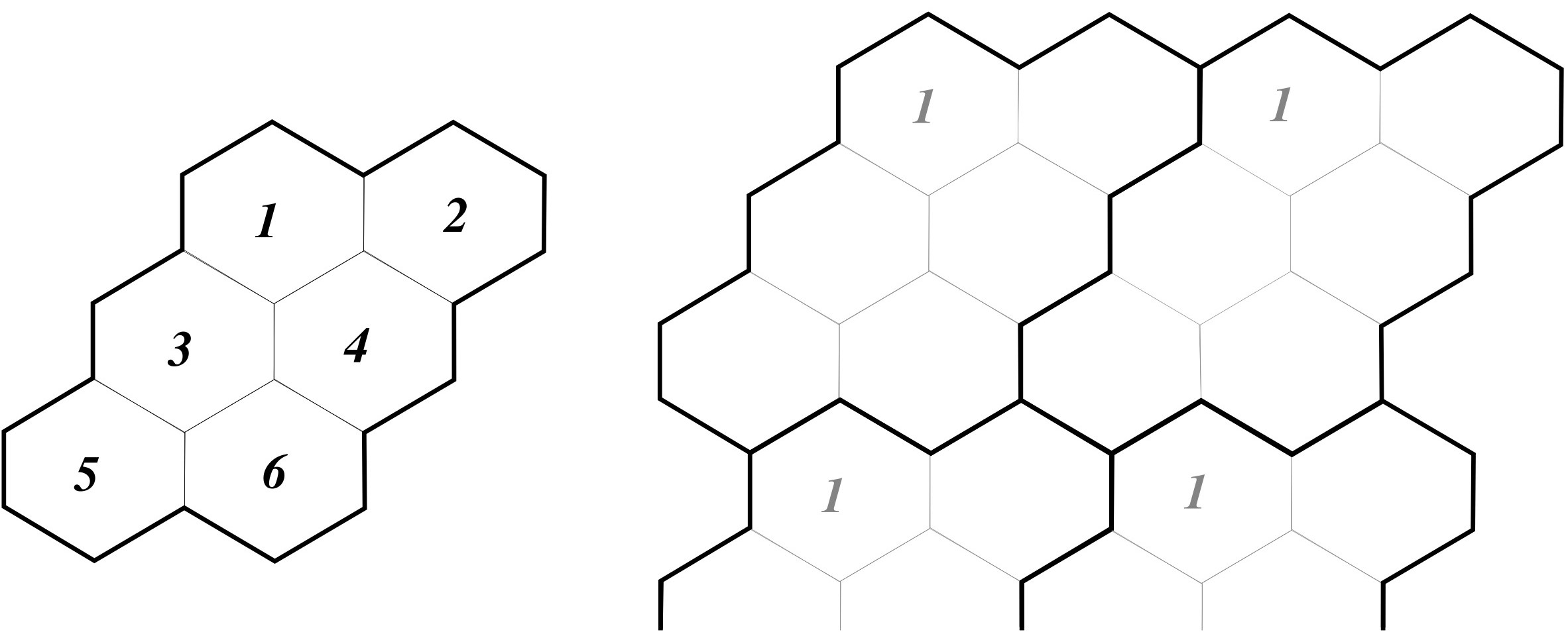}
\caption{\small{A tiling of the plane with translates of hexagon $H$ defines a periodic $6$-coloring}}
\label{color classes}
\end{figure}

To avoid unit $C$-distance in $H$, remove the boundary points lying between the points $A_1$ and $A_4$, including $A_4$ but not $A_1$.
In this way no antipodal point pair is monochromatic.
For simplicity call the resulting half-open hexagon still $H$.
Now consider a tiling of the plane by translates of $H$ and assign colors $1$ through $6$ periodically as shown in Figure \ref{color classes}. We can assume that the centers of the hexagons form a lattice, that we denote by $\mathcal{L}$.

Note that the main difference between this coloring scheme and Chilakamarri's general $7$-coloring is that here we can take advantage of $C$ not being strictly convex: in the direction perpendicular to the sides shared with $C/2$ the monochromatic hexagons can be placed such that they are separated by only one differently colored hexagon.\\

Now we justify that our coloring is proper:
as mentioned earlier a unit $C$-distance is not realized within the hexagons. All is left is to check that two hexagons of the same color are too far from each other to determine a point pair unit $C$-distance apart. As the color classes are congruent, it is enough to verify the statement for one specific color class, say the class of red points. We can also assume that one of the red hexagons has the origin as its center, thus the set of centers of red hexagons form a sublattice $\mathcal{L'}$.
By the symmetry of $C$ it is enough to show that polygons $ \mathcal{L'} + C/2 \oplus H$ form a packing, where $\oplus$ denotes the Minkowski sum of the two polygons, that is:
$$C/2 \oplus H= \{  c+h \text{ } | \text{ }  c \in C/2, \text{ } h \in H \}. $$

\begin{figure}[htp]
   \centering
    \includegraphics[width=12cm]{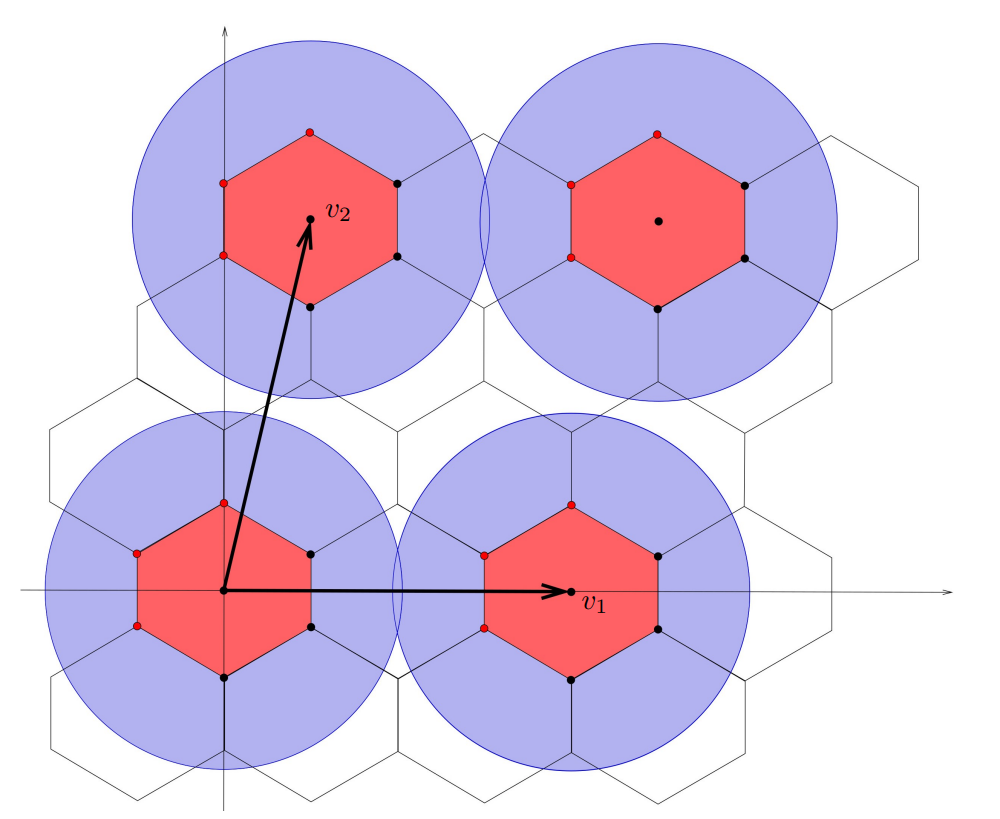}
    \caption{\small{ For any lattice vector $v$ of $\mathcal{L'} \setminus \{ \pm v_1, \underline{0} \}$ the Euclidean unit circle $B_2$ is disjoint from all translates $B_2+ v$}}
    \label{8_kor}
    \end{figure}

Straightforward calculations finish the proof.
Without loss of generality we can assume that $C$ has circumradius one.
Let $v_1$ and $v_2$ denote the basis vectors of the lattice $\mathcal{L'}$ where we can assume that $v_1$ is perpendicular to the sides shared with $C/2$. Then for any vector $\lambda \in \mathcal{L'} $, polygons $\lambda + C/2 \oplus H$ and $\lambda \pm v_1  + C/2 \oplus H$ are trivially disjoint.
From definition $H+C/2 \subseteq C$ so $H+C/2$ also has circumradius at most one. Hence it is enough to check that with the exception of $\underline{0}$ and $\pm v_1$ any lattice vector has Euclidean length at least two. This obviously holds (see Figure \ref{8_kor}),
thus the coloring is indeed proper.

\newpage
Now consider regular polygons with greater number of vertices. We show that almost the same coloring scheme works for all the remaining cases: Two sides of $H$ can always be two opposite sides of $C/2$, we only have to be careful with the choice of the remaining two vertices.
For the case of a regular $10$- and $12$-gons choosing the halving points on the boundary line of $C/2$ still works, we can simply define hexagon $H$ as shown in Figure \ref{hatszog_10,12}.

\vspace{-0.2cm}

\begin{figure}[htp]
\centering
\includegraphics[width=8.6cm]{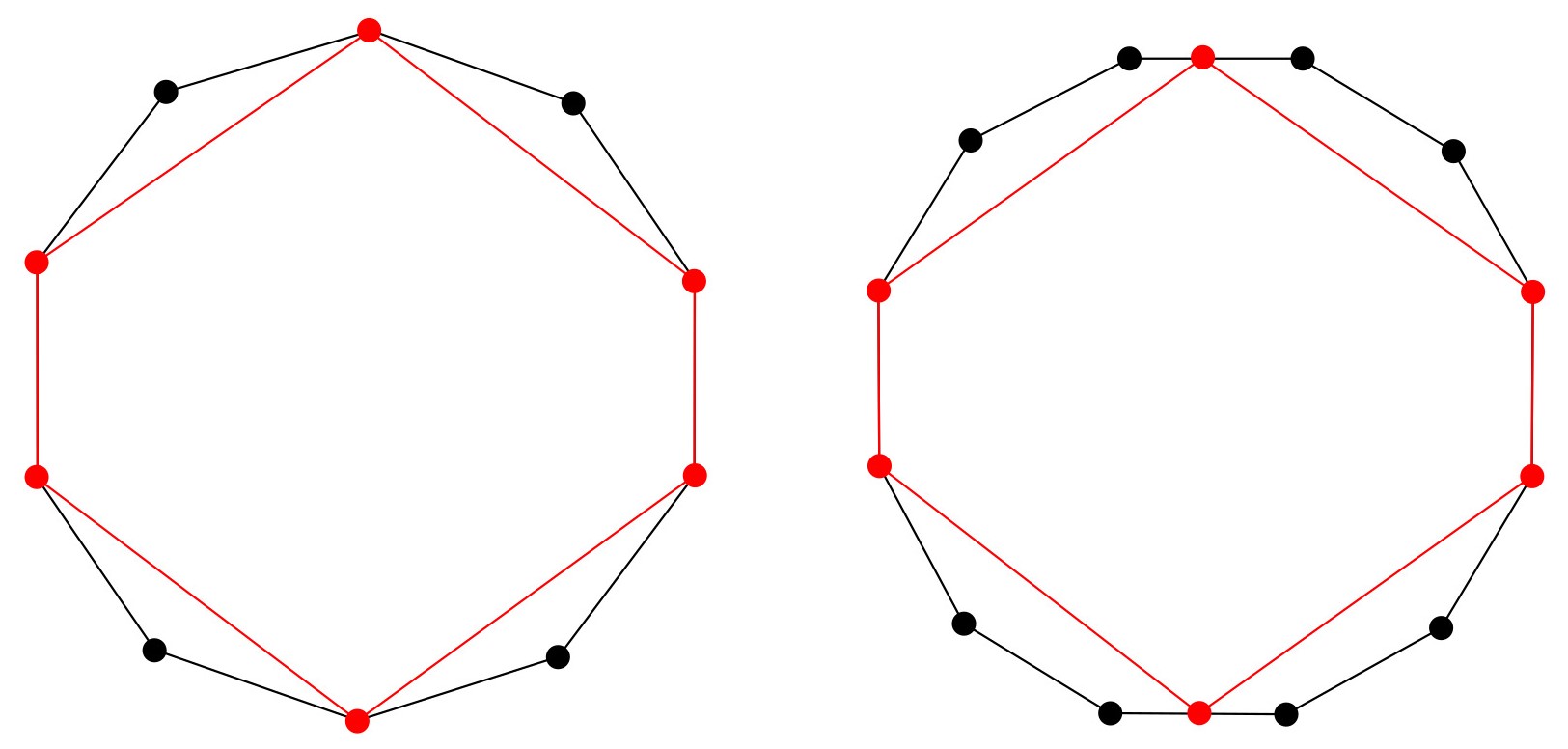}
\caption{\small{The centrally symmetric hexagons inscribed in $C/2$ chosen in the case of the regular $10$- and $12$-gon}}
\label{hatszog_10,12}
\end{figure}

\begin{figure}[htp]
\centering
\includegraphics[width=12.8cm]{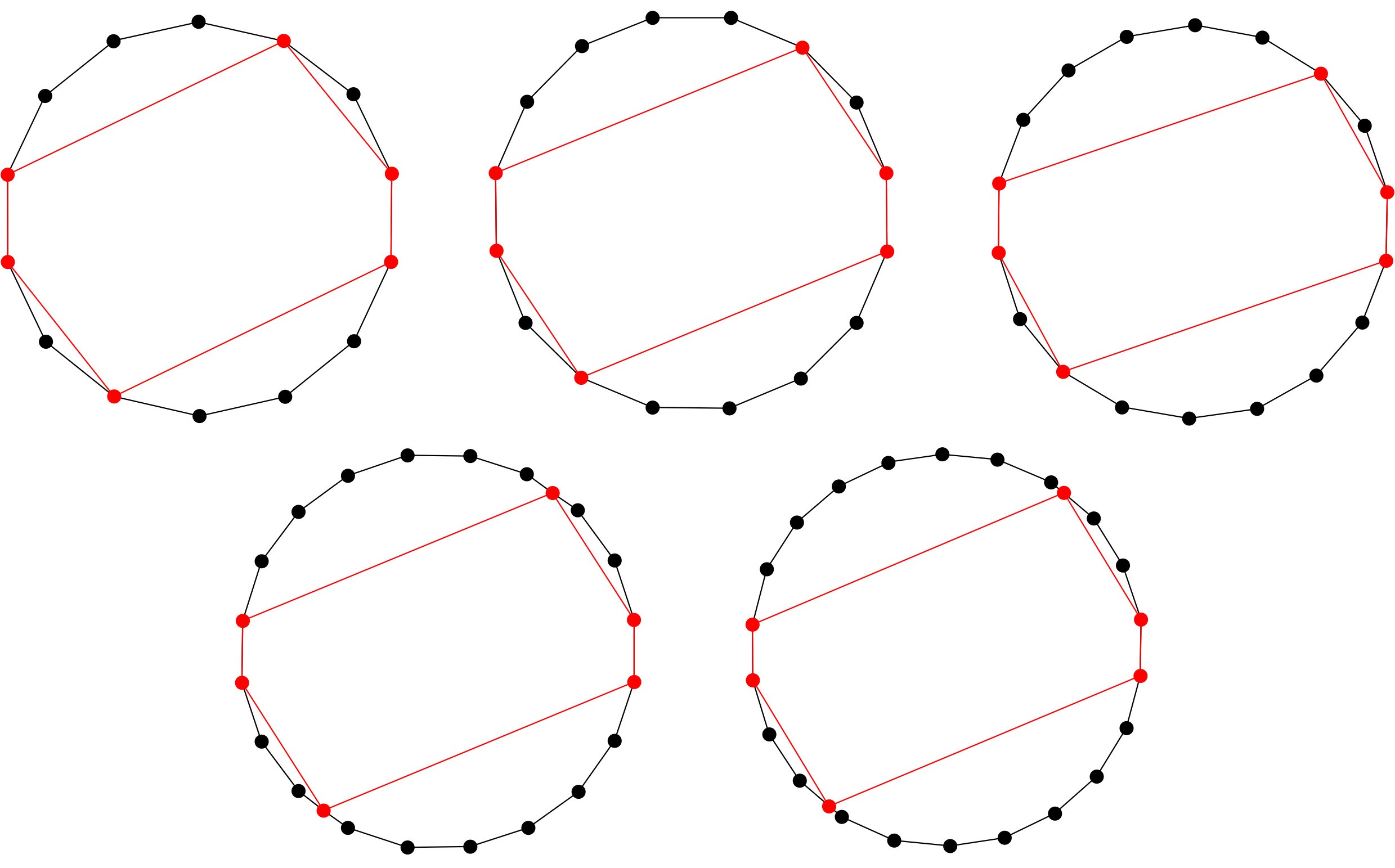}
\caption{\small{The centrally symmetric hexagon $H$ inscribed in $C/2$ chosen for the regular $14$-, $16$-, $18$-, $20$- and $22$-gon}}
\label{hatszogek_14-22}
\end{figure}

However, in the remaining cases we had to flatten the hexagons in order to get a proper coloring. In the case of regular $14$-, $16$- and $18$-gons some other vertices of $C/2$ were chosen. But in the final two cases only non-vertex points seemed to be working:
for $n=20$ bisectors of some other sides were chosen, and
for $n=22$ we divided one of the sides in the ratio $0.68:0.32$. For details see Figure \ref{hatszogek_14-22}.

Proving that the colorings defined by the above hexagons are proper needs more and more computations as the number of vertices grows. Since all cases are similar, we will only give the details of the calculations in the regular dodecagon's case, in Section \ref{proof} and in the regular $22$-gon's case in the Appendix.


\section{The regular dodecagon's case} \label{proof}
Now we present the details of the proof in the case of the regular dodecagon. Consider the regular dodecagon centered at the origin with circumradius $2$, whose vertices are:
$$\big( \pm 1, \pm \sqrt{3} \big), 
\big( \pm \sqrt{3}, \pm 1 \big), 
\big( \pm 2, 0 \big),
\big( 0, \pm 2 \big).$$

\vspace{0.1cm}

\begin{figure}[h]
\centering
\includegraphics[width=11.5cm]{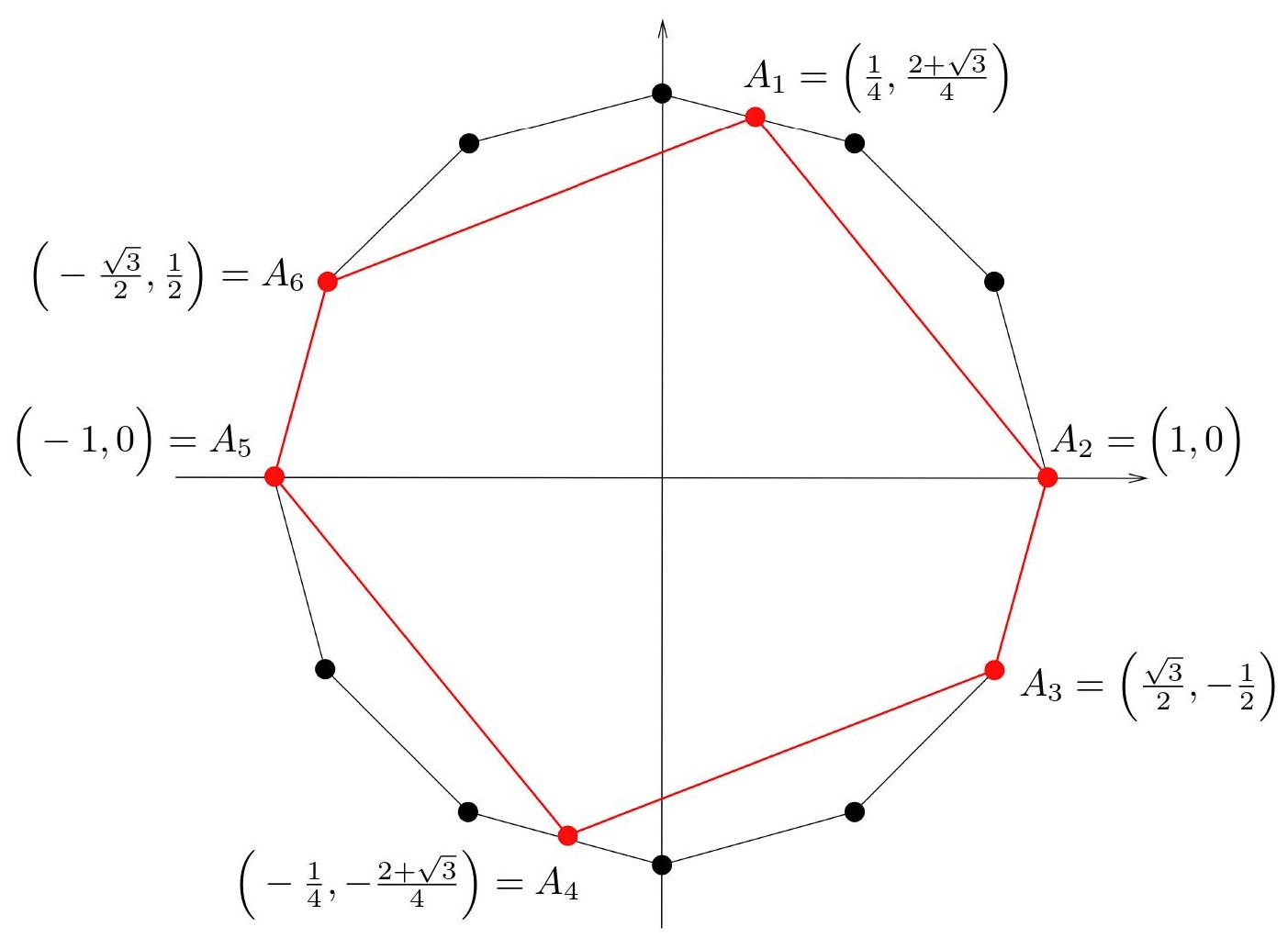}
\vspace{0.2cm}

\caption{\small{$H$ hexagon inscribed in $C/2$}}
\label{12_1}
\end{figure}

Let $H$ be the symmetric hexagon inscribed in $C/2$ as defined in Section \ref{coloring}: take two opposite sides of $C/2$, for example the sides parallel to vector $( 2- \sqrt{3}, 1  )$ and choose the two additional points such that they halve parts of the boundary line of $C/2$ between the chosen sides. Denote these six vertices by $A_i$ ($i = 1 \dots 6$) in a clockwise order as shown in Figure \ref{12_1}. As before, let $H$ be the half-open hexagon defined by the points $A_i$ that does not contain the line segment connecting the points $A_1$ and $A_4$ and the point $A_1$ itself.

Therefore the hexagonal tiling of the plane with hexagon $H$ is the packing by Voronoi regions of the lattice $\mathcal{L}$ spanned by vectors $\Big( \frac{1-2 \cdot \sqrt{3}}{4}, \frac{4+ \sqrt{3}}{4} \Big)$ and $\Big( \frac{2+\sqrt{3}}{2}, - \frac{1}{2} \Big)$.
The basis vectors of the sublattice $\mathcal{L'}$ corresponding to the single color class containing the hexagon centered at the origin are:
\begin{itemize}
    \item  $v_1= \big( \frac{3- 6 \cdot \sqrt{3}}{4}, \frac{12 +3 \cdot \sqrt{3}}{4} \big)$,
    \item $v_2= \big( 2+ \sqrt{3}, -1 \big)$.
\end{itemize}

\vspace{0.2cm}

As mentioned in Section \ref{coloring} what we need to show is that polygons $\mathcal{L'} + C/2 \oplus H$ form a packing.

The vertices of $C/2 \oplus H$ are:
$$B_1= \Bigg( \frac{1}{4}, \frac{6+ \sqrt{3}}{4} \Bigg), \text{ }
B_2= \Bigg( \frac{3}{4}, \frac{2+3 \cdot \sqrt{3}}{4} \Bigg), \text{ }
B_3= \Bigg( \frac{1+ 2 \cdot \sqrt{3}}{4}, \frac{4+ \sqrt{3}}{4} \Bigg), \text{ }
B_4= \Bigg( \frac{2+ \sqrt{3}}{2}, \frac{1}{2} \Bigg),$$
$$ B_5= \Big( 2, 0 \Big), \text{ }
B_6= \Bigg( \sqrt{3}, -1\Bigg),  \text{ }
B_7= \Bigg(  \frac{1+\sqrt{3}}{2}, - \frac{1+\sqrt{3}}{2} \Bigg), \text{ }
B_8= \Bigg( \frac{1}{4}, - \frac{2+3 \cdot \sqrt{3}}{4} \Bigg) \text{ etc.}$$

\begin{figure}[htp!]
\centering
\hspace{-1cm}
\includegraphics[width=12.6cm]{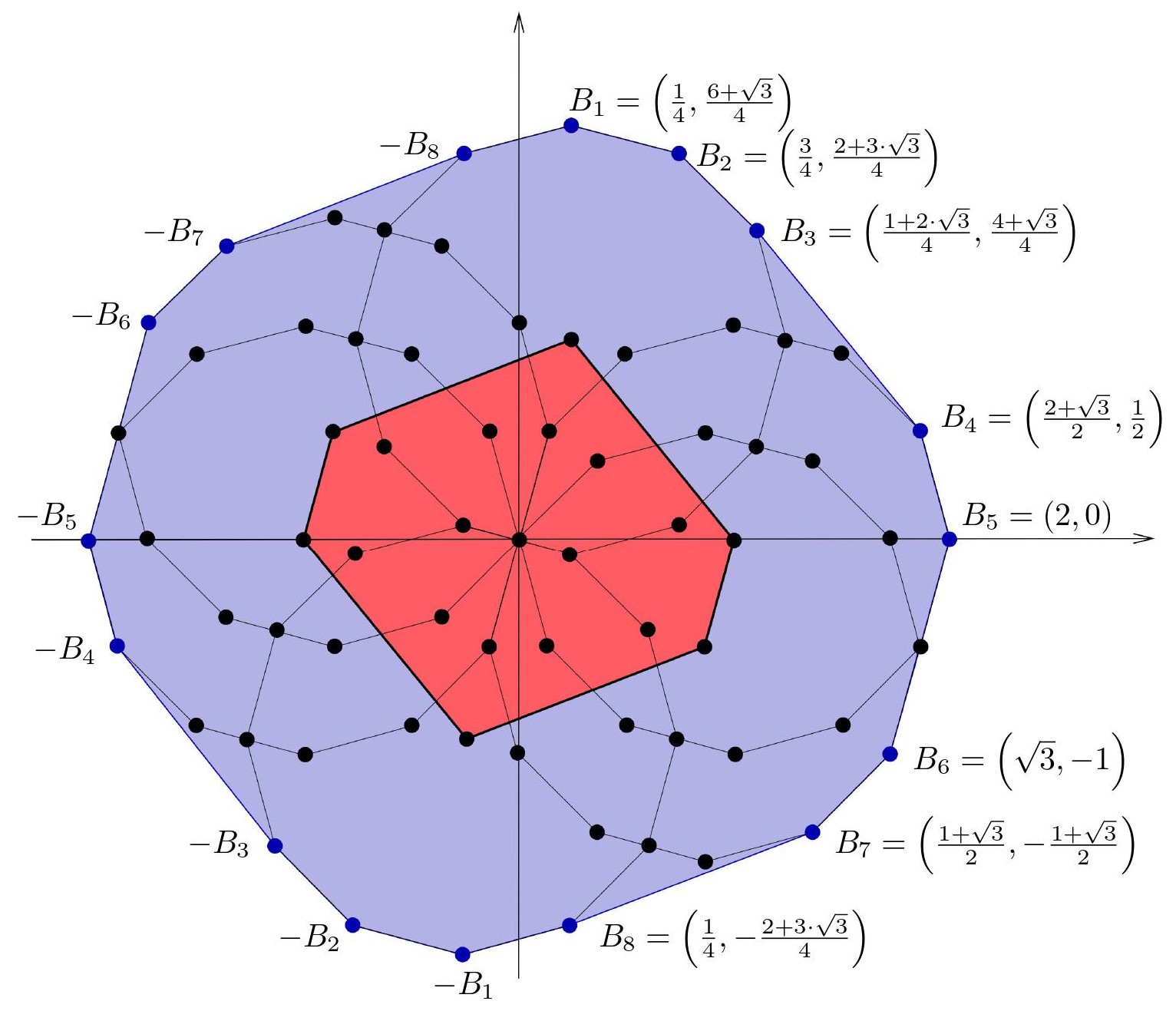}

\vspace{0.1cm}
\caption{\small{Minkowski sum of $H$ and $C/2$}}
\label{Minkowski}
\end{figure}

The coordinates of the remaining vertices can be obtained by symmetry -- see Figure \ref{Minkowski}.

As the coloring is regular, it is enough to pick hexagon $H$ centered at the origin and show that $H\oplus1/2 C$ is disjoint from $\lambda' + H \oplus 1/2 C$ for all $\lambda' \not\equiv 0$ in $\mathcal{L'}$.

By definition $H \oplus C/2$ has circumradius $2$.
Inside the circle of radius $4$ centered at the origin there are $4$ lattice points of $\mathcal{L'}$ besides the origin and by symmetry we only have to check $2$ of them and the corresponding hexagons, namely:
\begin{itemize}
    \item $H_1:= H+v_2$ and
    \item $H_2:=H+v_1+v_2$.
    
\end{itemize}

$H$ and $H_1$ are separated by exactly one differently colored hexagon which is enough as $v_2$ is perpendicular to the common sides of $H$ and $C/2$.
All is left is to give a line that separates $H$ from $H_2$.
For example consider the line $l$ defined by equation:

$$ y = - \frac{2+ \sqrt{3}}{3}x + \frac{5 + 2 \cdot \sqrt{3}}{3}.$$

\vspace{0.2cm}

   \begin{figure}[!htp]
    \centering
    \includegraphics[width=16.7cm]{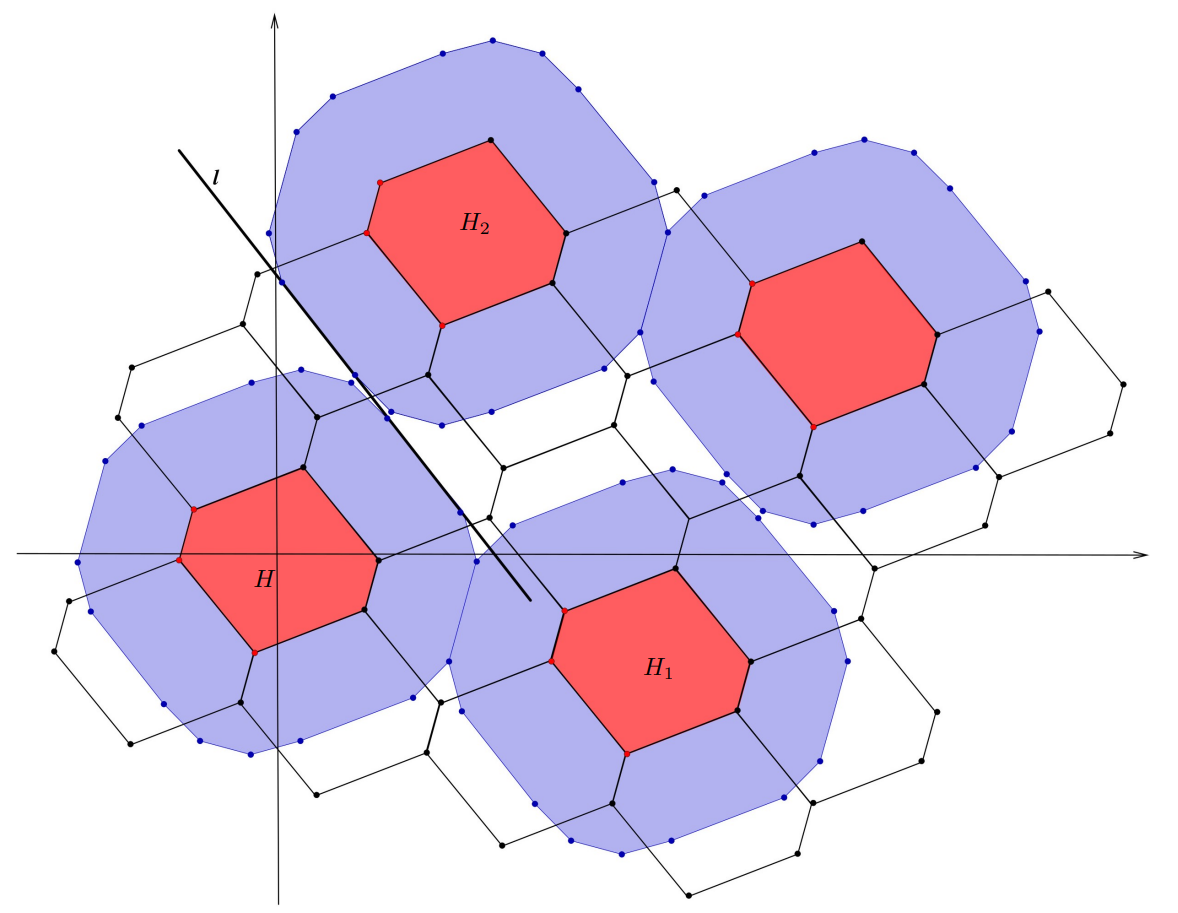}
    \caption{\small{Line $l$ separates $H \oplus C/2$ and $H_2 \oplus C/2$}}
    \label{12_szog}
    \end{figure}

It is straightforward to check that $l$ goes through two parallel sides of $H \oplus C/2$ and $H_2 \oplus C/2$ (which are on the opposite sides of their centers) and the remaining vertex points of $H \oplus C/2$ are below line $l$, while all of the remaining vertex points of $H_2 \oplus C/2$ are above it (see Figure \ref{12_szog}). Therefore  $H \oplus C/2$ and $H_2 \oplus C/2$ are disjoint, thus the coloring is proper.\\

\begin{figure}[!htp]
    \centering
    \includegraphics[width=15cm]{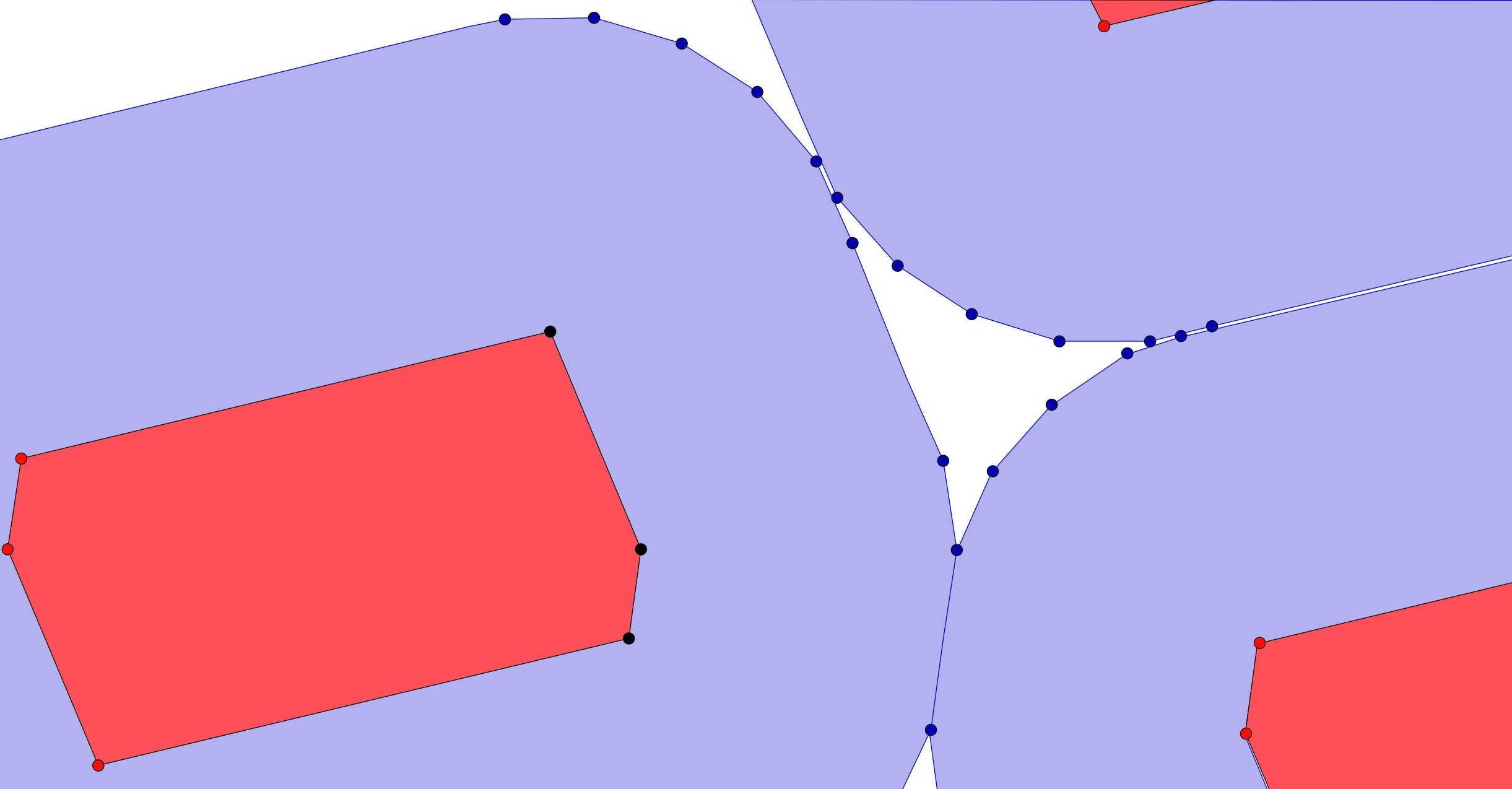}
    \caption{\small{In the $6$-coloring of the Minkowski plane equipped with the regular $22$-gon metric\\
    monochromatic hexagons get dangerously close together}}
    \label{22}
    \end{figure}

We remark that in the presented example one can define hexagon $H$ in many different ways as the coloring scheme is quite flexible in this case. However, as the number of vertices increases, the range of possible choices narrows down quickly. For example in the case of the regular $22$-gon we have to be really carefull with the definition of hexagon $H$: as Figure \ref{22} shows, our coloring is almost rigid.

\section{An asymmetric Ramsey-type problem} \label{Szlam}

Another direction for generalizing the Hadwiger-Nelson problem is to replace the pair of points at unit distance by another finite point configuration. Moreover we can look for different configurations in each color class.
Solving a problem raised in \cite{Erdos2} Juhász \cite{Juhasz} showed that in any red-blue coloring of the plane there are either two red points distance one apart or there is a blue congruent copy of any configuration with at most $4$ points. In the opposite direction Juhász also showed that her theorem does not remain true if we replace $4$ by $12$. Csizmadia and Tóth \cite{Csizmadia-Toth} later improved the above result, they proved that $4$ can not even be increased to $8$.\\

In the rest of the paper we are interested in the following variation of the question above: is it true for a given point configuration $K\subseteq \mathbb{R}^n$ that in any red-blue coloring of the $n$-dimensional Euclidean space there are either two red points distance one apart or there is a blue translate of configuration $K$? Let $k_n$ denote the largest value $k$ such that the answer is `yes' for all configuration $K$ of at most $k$ points. Note that the $2$-dimensional case is not the same as it was in the original problem since `congruent copy' has been replaced by \textit{`translate'}.
\\
This variant of the problem was first considered in \cite{Szlam} by Szlam who showed that $k_n$ grows exponentially in the number of dimensions. More precisely, Szlam's theorem states that $\chi(\mathbb{R}^n)$, the chromatic number of the $n$-dimensional Euclidean space (i. e. the smallest number of colors that are needed to color $\mathbb{R}^n$ so that no two points of the same color determine unit distance) provides a lower bound on the value of $k_n$.
For the sake of completeness, we include his short proof as well.

\begin{lemma} [Szlam I. \cite{Szlam}] \label{SzlamI}
Assume that there exists a $k$-point configuration $K$ and a red-blue coloring of $\mathbb{R}^n$ such that the red color class avoids unit distance and the blue color class avoids all translates of $K$. Then $\mathbb{R}^n$ can be properly $k$-colored that is none of the $k$ color classes contains unit distance. Hence $k_n +1 \geq \chi(\mathbb{R}^n)$.
\begin{proof}
Assume that we are given a configuration $K= \{a_1, \dots a_k \}$ and a red-blue coloring of $\mathbb{R}^n$ with the desired properties.
Then, for each $x \in \mathbb{R}^n$ there is at least one index $i$ such that $x + a_i$ is red. Now let us color the point $x$ with color number $i$: as there are no red points unit distance apart, this indeed defines a proper $k$-coloring.
\end{proof}
\end{lemma}

The famous theorem of Frankl and Wilson \cite{Frankl-Wilson} gives an exponential lower bound for the chromatic number of the Euclidean space: it states that in case $n$ tends to infinity
$$ \chi(\mathbb{R}^n )\geq (1.2+o(1))^n.$$
Hence, by Lemma \ref{SzlamI} if $n$ tends to infinity, every red-blue coloring of $\mathbb{R}^n$  contains a blue translate of all $k$ point configuration when $k < (1.2+o(1))^n$.\\

Szlam also a gave a partial converse to Lemma \ref{SzlamI} that considers only regular colorings:

\begin{lemma} [Szlam II. \cite{Szlam}] \label{SzlamII}
Assume that $\mathbb{R}^n$ can be properly $k$-colored by a regular coloring, with color classes $C_i=C_1+v_i$ (for $i=1 \dots k$).
Then there exists a red-blue coloring of $\mathbb{R}^n$ and a $k$-point configuration, namely $K= \{ v_1, v_2, \dots v_k \}$ such that the red color class avoids unit distance and the blue color class avoids all translates of $K$.
\end{lemma}

Since Szlam's original paper was published many applications and generalizations of his work were considered (see for example \cite{Burkert-Johnson, Johnson-Szlam, Krizan, Weatherspoon}). Here we need the analogous result considering Minkowski spaces:
Let $k_n(C)$ denote the largest $k$ value such that in any red-blue coloring of the Minkowski space determined by $C$ there are either two red points distance one apart or there is a blue translate of any configuration with at most $k$ points. An easy observation is that both Lemma \ref{SzlamI} and Lemma \ref{SzlamII} can be extended to Minkowski spaces (as noted e.g. in \cite{Burkert-Johnson}). As the $6$-colorings described in Theorem \ref{theorem1} are a regular, we immediately get the following corollary:

\begin{corollary}
Let $(\mathbb{R}^2, C)$ be a Minkowski plane whose unit circle is a regular polygon with an even number of vertices. In case $C$ has at most $22$ vertices
then there exists a red-blue coloring of the plane and a configuration $K$ of $6$ points such that there is no red point pair unit $C$-distance apart, and the blue color class avoids all translates of $K$.
\end{corollary}

We finish the paper with a small remark on Szlam's results: We noticed that although the proof of Lemma \ref{SzlamII} is short and straightforward, it can be a bit misleading. To see its inconvenience notice how the proper coloring defined in Lemma \ref{SzlamI} is not regular: color classes are generated by covering the plane with translates of the unit distance avoiding red set. Call a coloring with such structure subregular. More precisely we call a proper $k$-coloring with color classes $C_1$ $\dots$ $C_k$ \textit{subregular} if there exist vectors $v_1$ $\dots$ $v_k$ such that $C_i$ is a subset of $C_1+v_i$.
We show that Lemma \ref{SzlamII} can be extended to subregular colorings in a very natural way:

\begin{theorem} \label{theorem2}
Let $(\mathbb{R}^n, C)$ be an $n$-dimensional Minkowski space. Assume $\mathbb{R}^n$ can properly be $k$-colored by a subregular coloring defined by a $C$-unit distance avoiding set $C_1$ and vectors $v_1$, $v_2$ $\dots$ $v_k$. Then there exists a red-blue coloring of $\mathbb{R}^n$ and a $k$-point configuration, namely $K= \{ -v_1, -v_2, \dots -v_k \}$ such that the red color class avoids unit $C$-distance and the blue color class avoids all translates of $K$.

\begin{proof}
Let the points of $C_1$ be colored red, and color all the remaining points blue. As promised, let us consider the configuration $K=\{ - v_1, -v_2, \dots, -v_k \}$. We wish to show that for an arbitrary vector $m$ color class $C_1$ contains at least one point of $K+m$.
Without loss of generality we can assume $v_1 \equiv 0$. Hence if $m \in C_1$ there is nothing to prove. Assume that $m \notin C_1$. In this case there exists an index $i$ such that $m \in C_1+v_i$ which leads to $-v_i+m \in C_1$.
\end{proof}
\end{theorem}
It follows that for all $n$ and $C$ the value $k_n(C)$ is exactly the smallest number $k$ such that there exists a subregular $k$-coloring of $(\mathbb{R}^n, C)$.
As we have seen, known proper colorings of Minkowski planes are usually regular. In the $3$-dimensional Euclidean a proper $15$-coloring was defined in \cite{R-Toth} and independently in \cite{Coulson} which give the best upper bound for the chromatic number of $\mathbb{R}^3$. Although these colorings are not rigorously regular, they can be turned into regular colorings in a trivial way. However, in higher dimensions proper colorings are typically only subregular. The best known upper bound on the chromatic number of the Euclidean space was established by Larman and Rogers \cite{Larman-Rogers} who identified a subregular proper coloring of $\mathbb{R}^n$ with $(3+o(1))^n$ color classes, meaning $k_n < (3+o(1))^n$. By Lemma \ref{SzlamII} this result is tight up to the constant factor. For Minkowski spaces the analogous theorem with $(4+o(1))^n$ color classes was proved by Kupavskii \cite{Kupavskii}, meaning $k_n(C) < (4+o(1))^n$.

\section*{Acknowledgement}
I would like to express my gratitude to Géza Tóth for his great help and for his many helpful remarks. I am also thankful to Dan Ismailescu for informing us about the current state of the problem.


\newpage
\section*{Appendix}
For the sake of completeness, we give the calculations for the regular $22$-gon's case: Consider the regular $22$-gon centered at the origin with circumradius $2$ so that the coordinates of the vertices of $C/2$ are:
$$\left(\pm \cos(\frac{2t\pi}{11}), \pm \sin(\frac{2t\pi}{11}) \right)$$
where $t=0, 1, 2, 3, 4, 5$.\\

Let $H$ be the symmetric hexagon inscribed in $C/2$ as defined in Section \ref{coloring}: take the sides of $C/2$ parallel to vector $\left( 1+ \cos(\frac{20 \pi}{22}), \sin(\frac{20 \pi}{22})  \right)$ and choose the two additional points such that they divide the sides parallel to vector $\left( \cos(\frac{6 \pi}{22})- \cos(\frac{4 \pi}{22}), \sin(\frac{6 \pi}{22})- \sin(\frac{4 \pi}{22}) \right)$ in the ratio $0.68:0.32$.\\
Coordinates of the vertices of $H$ are:
\begin{itemize}
    \item
$A_1=-A_4= \left( 0.32 \cos(\frac{4 \pi}{22}) + 0.68 \cos(\frac{6 \pi}{22}), 0.32 \sin(\frac{4 \pi}{22}) + 0.68 \sin(\frac{6 \pi}{22}) \right)$,
\item
$A_2=-A_5= \left(1, 0 \right)$,
\item
$A_3=-A_6= \left( \cos(\frac{2 \pi}{22}), \sin(\frac{-2 \pi}{22}) \right)$.
\end{itemize}

As before, let $H$ be the half-open hexagon defined by the points $A_i$ that does not contain the line segment connecting the points $A_1$ and $A_4$ and the point $A_1$ itself.
For $i=1 \dots 6$ let $a^i$ denote the position vector of point $A_i$. Then the hexagonal tiling of the plane with hexagon $H$ is the packing by Voronoi regions of the lattice $\mathcal{L}$ spanned by vectors $ \left(a^1 + a^6 \right)$ and $ \left( a^2 + a^3\right)$. The basis vectors of the sublattice $\mathcal{L'}$ corresponding to the single color class containing the hexagon centered at the origin are:
\begin{itemize}
    \item  $v_1= \left( 3(a^1+a^6) \right)$,
    \item $v_2= \left( 2 (a^2 +a^3)  \right)$.
\end{itemize}

Once again, we need to show is that polygons $\mathcal{L'} + C/2 \oplus H$ form a packing. The vertices of $C/2 \oplus H$ are:

  \begin{multicols}{2}
\begin{itemize}

    \item $\pm \left( \cos(\frac{12 \pi}{22}) + a^1_x, \text{ } \sin(\frac{12 \pi}{22}) + a^1_y \right)$

    \item $ \pm \left( \cos(\frac{10 \pi}{22}) + a^1_x, \text{ } \sin(\frac{10 \pi}{22}) + a^1_y \right) $

    \item $ \pm \left( \cos(\frac{8 \pi}{22}) + a^1_x, \text{ } \sin(\frac{8 \pi}{22}) + a^1_y \right) $

    \item $ \pm \left( \cos(\frac{6 \pi}{22}) + a^1_x, \text{ } \sin(\frac{6 \pi}{22}) + a^1_y \right)$

    \item $ \pm \left( \cos(\frac{4 \pi}{22}) + a^1_x, \text{ } \sin(\frac{4 \pi}{22}) + a^1_y \right) $

    \item $ \pm \left( \cos(\frac{2 \pi}{22}) +a^1_x, \text{ } \sin(\frac{2 \pi}{22}) + a^1_y\right)$

    \item $ \pm \left( \cos(\frac{2 \pi}{22}) + a^2_x, \text{ } \sin(\frac{2 \pi}{22})+ a^2_y \right) $

    \item $ \pm \left( 1 + a^2_x, \text{ } a^2_y \right) =(2,0) $

    \item $ \pm \left( \cos(\frac{2 \pi}{22})+ a^3_x, \text{ } \sin(\frac{-2 \pi}{22})+ a^3_y \right)$

    \item $ \pm \left( \cos(\frac{4 \pi}{22})+ a^3_x, \text{ } \sin(\frac{-4 \pi}{22}) + a^3_y \right) $

    \item $ \pm \left( \cos(\frac{6 \pi}{22}) + a^3_x, \text{ } \sin(\frac{-6 \pi}{22}) + a^3_y \right)$

    \item $ \pm \left( \cos(\frac{8 \pi}{22}) + a^3_x, \text{ } \sin(\frac{-8 \pi}{22}) + a^3_y \right)$

    \item $ \pm \left( \cos(\frac{10 \pi}{22}) + a^3_x, \text{ } \sin(\frac{-10 \pi}{22})+ a^3_y \right) $

\end{itemize}
  \end{multicols}

Since  $H \oplus C/2$ has circumradius $2$, it is enough to consider lattice points of $\mathcal{L'}$  with length at most 4. Besides the origin there are $6$ such points and by symmetry we only have to check $3$ of them and the corresponding hexagons, namely:\\
 $H_1:= H+v_2$,
 $H_2:=H+v_1$ and $H_3:=H+v_1+v_2$.

$H$ and $H_1$ are separated by exactly one differently colored hexagon which is enough as $v_2$ is perpendicular to the common sides of $H$ and $C/2$.
All is left is to give a line $l_1$ that separates $H$ from $H_2$ and a line $l_2$ that separates $H$ from $H_3$.
For example consider the lines $l_1$ and $l_2$ defined by the following equations:

\begin{align*}
    y- \sin(\frac{12 \pi}{22}) -  a^6_x - \frac{1}{300} &= \frac{ a^6_y - a^1_y }{a^6_x - a^1_x} \left(x -\cos(\frac{12 \pi}{22}) - a^6_x \right)\\
    y- \sin(\frac{4 \pi}{22})- a^1_y - \frac{1}{50} &= \frac{\sin(\frac{4 \pi}{22})- \sin(\frac{2 \pi}{22})}{\cos(\frac{4 \pi}{22})- \cos(\frac{2 \pi}{22})}  \left(  x- \cos(\frac{4 \pi}{22}) - a^1_x \right)
\end{align*}

    \begin{figure}[!htp]
    \centering
    \includegraphics[width=17cm]{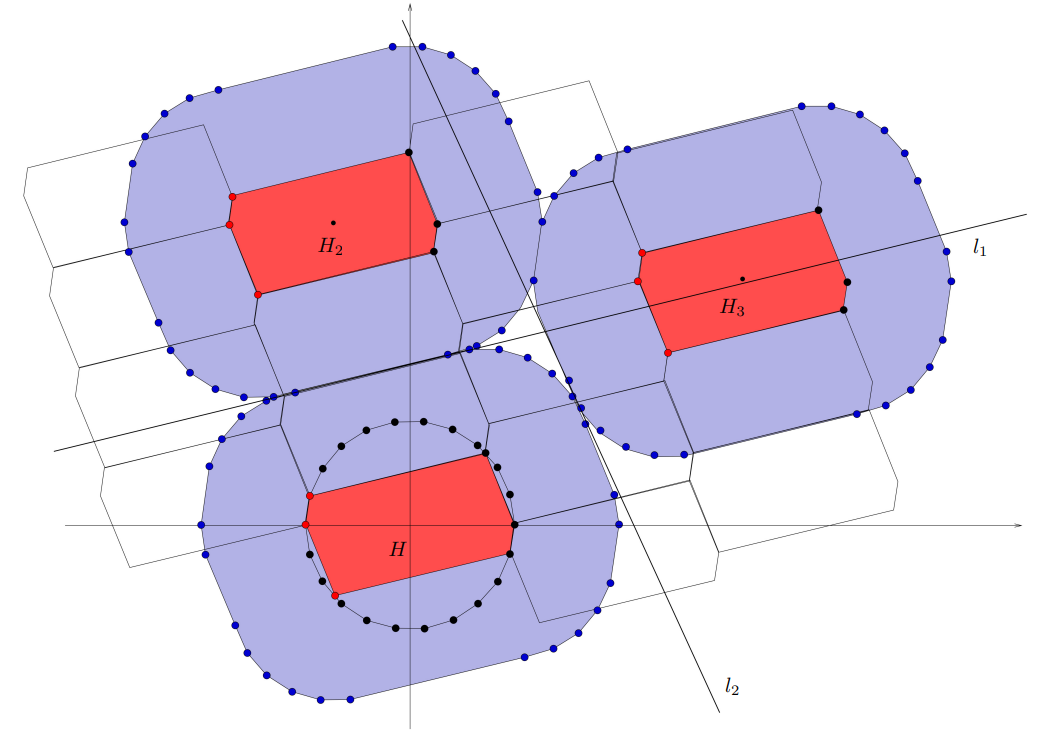}
    \caption{Lines $l_1$ and $l_2$ separate $H$ from $H_2 \oplus C/2$ and from $H_3 \oplus C/2 $}
    \label{22_2}
    \end{figure}

It is straightforward to check that all vertex points of $H \oplus C/2$ are below line $l_1$ and line $l_2$, while all vertex points of $H_2 \oplus C/2$ are above $l_1$, and all vertex points of $H_3 \oplus C/2$ are above $l_2$ (see Figure \ref{22_2}). Thus the coloring is indeed proper.
\end{document}